\documentclass[12pt]{article}

\usepackage[centertags]{amsmath}
\usepackage{amsfonts}
\usepackage{amssymb}
\usepackage{amsfonts}
\usepackage{latexsym}
\usepackage{amsthm}

\newcounter{abcd}
\Alph{abcd}

\newtheorem*{pred0}{Theorem}

\newcommand {\Hol}{\mathop{\rm Hol}\nolimits}

\newtheorem{prop}{Proposition}

\newtheorem*{corol}{Corollary}
\newcommand{\pr}{\noindent{\bf Proof.}\quad }

\newcommand{\epr}{\ $\blacksquare$  \vspace{2mm} }

\renewcommand{\Re}{\mathop{\rm Re}\nolimits}

\begin{document}

\begin{titlepage}

\begin{center}


{\bf\Large A rigidity theorem for holomorphic generators on the
Hilbert ball}

\bigskip

{\large Mark Elin\\}  {\footnotesize Department of Mathematics
\\ ORT  Braude College
\\ P.O. Box 78, 21982 Karmiel, Israel
\\ e-mail: mark.elin@gmail.com}

\bigskip

{\large Marina Levenshtein\\}  {\footnotesize Department of
Mathematics
\\ The Technion --- Israel Institute of Technology
\\ 32000 Haifa, Israel
\\ e-mail: marlev@list.ru}

\bigskip

{\large Simeon Reich\\}  {\footnotesize Department of Mathematics
\\ The Technion --- Israel Institute of Technology
\\ 32000 Haifa, Israel
\\ e-mail: sreich@tx.technion.ac.il}

\bigskip

{\large David Shoikhet\\}  {\footnotesize Department of
Mathematics
\\ ORT  Braude College
\\ P.O. Box 78, 21982 Karmiel, Israel
\\ e-mail: davs27@netvision.net.il}

\end{center}

\bigskip

{\it \noindent 2000 Mathematics Subject Classification: 30C45,
30D05, 46T25, 47H20. }

\end{titlepage}

\begin{abstract}
We present a rigidity property of holomorphic generators on the
open unit ball $\mathbb{B}$ of a Hilbert space $H$. Namely, if
$f\in\Hol (\mathbb{B},H)$ is the generator of a one-parameter
continuous semigroup $\left\{F_t\right\}_{t\geq 0}$ on
$\mathbb{B}$ such that for some boundary point
$\tau\in\partial\mathbb{B}$, the admissible limit
$K$-$\lim\limits_{z\rightarrow\tau}\frac{f(x)}{\|x-\tau\|^{3}}=0$,
then $f$ vanishes identically on $\mathbb{B}$.
\end{abstract}

Let $H$ be a complex Hilbert space with inner product
$\langle\cdot ,\cdot\rangle$ and induced norm $\|\cdot\|$. If $H$
is finite dimensional, we will identify $H$ with $\mathbb{C}^{n}$.
We denote by $\Hol(D,E)$ the set of all holomorphic mappings on a
domain $D\subset H$ which map $D$ into a subset $E$ of $H$, and put
$\Hol(D):=\Hol(D,D)$.

We are concerned with the problem of finding conditions for a mapping 
$F \in \Hol (D,E)$ to coincide identically with a given holomorphic
mapping on $D$, when they behave similarly in a neighborhood of a
boundary point $\tau\in\partial D$.

For holomorphic self-mappings of the open unit disk ${\Delta
:=\{z\in\mathbb{C} :  |z|<1 \}}$, the following result in this
direction is due to  D. M. Burns and S. G. Krantz~\cite{B-K}.

\begin{prop}
Let $F\in\Hol (\Delta)$. If the unrestricted limit
$$\lim\limits_{z\rightarrow 1}\frac{F(z)-z}{(z-\tau)^{3}}=0$$ for
some $\tau\in \partial\Delta$, then $F\equiv I$ on $\Delta$.
\end{prop}

This assertion also holds when the unrestricted limit is replaced
with the angular one (see \cite{T-V} and \cite{B-T-V}). Recall
that a function $f\in\Hol (\Delta ,\mathbb{C})$ has an angular
limit $L:=\angle\lim\limits_{z\rightarrow\tau} f(z)$ at a point
$\tau\in\partial\Delta \, $ if $\, f(z)\rightarrow L$ as
$z\rightarrow\tau$ in each nontangential approach region $$\Gamma
_{k}(\tau):=\left\{z\in\Delta : \frac{|z-\tau|}{1-|z|}<k\right\},
\, \, k>1. $$ In this case it is convenient to set
$f(\tau):=\angle\lim\limits_{z\rightarrow\tau} f(z)$. Moreover, in
a similar way, one defines the angular derivative of $f$ at
$\tau\in\partial\Delta$ by
$f'(\tau):=\angle\lim\limits_{z\rightarrow\tau}
\frac{f(z)-f(\tau)}{z-\tau}$.

A point $\tau\in\overline{\Delta}$ is a fixed point of
$F\in\Hol(\Delta)$ if either $F(\tau)=\tau$, where $
\tau\in\Delta$, or $\lim\limits_{r\rightarrow 1^-}F(r\tau)=\tau$,
where $\tau\in\partial\Delta=\{z:|z|=1\}$. If $F$ is not an
automorphism of $\Delta$ with an interior fixed point, then by the
classical Schwarz--Pick lemma and the Julia--Wolff--Carath\'eodory
theorem, there is a unique fixed point $\tau\in\overline{\Delta}$
such that for each $z\in\Delta$,
$\lim\limits_{n\rightarrow\infty}F_{n}(z)=\tau$, where the $n$-th
iteration $F_{n}$ of $F$ is defined by $F_1=F, \, F_{n}=F\circ
F_{n-1},\ n=2,3,\ldots$. This point is called the Denjoy--Wolff
point of $F$.  Moreover, a boundary fixed point
$\tau\in\partial\Delta$ of $F$ is its Denjoy--Wolff point if and
only if $F'(\tau)\in(0,1]$.

A rigidity result for generators of
one-parameter continuous semigroups on $\Delta$ (see Proposition 2
below) has been proved in \cite{E-L-S-T} and \cite{L-R-S}. To
formulate it, we first recall the definitions of these notions.

Let $D\subset H$ be a domain in the Hilbert space $H$. We say that
a family $S=\left\{F_t\right\}_{t\geq 0}\subset\Hol (D)$ is {\bf a
one-parameter continuous semigroup on $D$} (a semigroup, in short)
if

(i) $F_{t}(F_{s}(z))=F_{t+s}(z)$ for all $t,s\geq 0$ and all $z \in D$,

\noindent and

(ii) $\lim\limits_{t\rightarrow 0^+}F_{t}(z)=z$ for all $z\in D$.

A semigroup $S=\left\{F_t\right\}_{t\geq 0}\subset\Hol (D)$ is
said to be generated if for each $z\in D$, there exists the strong
limit  $$ \lim_{t\rightarrow 0^{+}}\frac{1}{t}[z-F_{t}(z)]=f(z).$$
In this case the mapping $f:D\mapsto H$ is called the {\bf
(infinitesimal) generator} of $S$.

A well-known representation of generators on $\Delta$ is due to E.
Berkson and H. Porta \cite{B-P}, namely:

\emph {A function $f\in\Hol (\Delta ,\mathbb{C})$ is a generator
if and only if there are a point $\tau\in\overline{\Delta}$ and a
function $p\in\Hol (\Delta , \mathbb{C})$ with $\Re p(z)\geq 0$
for all $z\in\Delta$ such that}
\begin{equation}\label{b11}
f(z)=(z-\tau)(1-\overline{\tau}z)p(z), \quad z\in\Delta.
\end{equation}
\emph {This point $\tau$ is the common Denjoy--Wolff point of the
semigroup generated by $f$. }

\begin{prop}
Let $g\in\Hol(\Delta , \mathbb{C})$ be the generator of a
one-parameter continuous semigroup. Suppose that
\begin{equation}\label{r1}
\angle\lim\limits_{z\rightarrow 1}\frac{g(z)}{|z-1|^{3}}=0.
\end{equation}
Then $g\equiv 0$ in $\Delta$.
\end{prop}

We take this opportunity to present here a completely different
proof of this assertion.

\pr
Suppose that $g$ does not vanish identically on $\Delta$.
Condition (\ref{r1}) implies that $\tau=1$ is the Denjoy--Wolff
point of the semigroup generated by $g$ (see Lemma 3 in
\cite{E-S}). So, $g$ has no null point in $\Delta$ (see Theorem 1
in \cite{E-S}). Consequently, $g$ can be represented by the
Berkson--Porta formula  $$g(z)=-(1-z)^{2}p(z), \quad z\in\Delta,$$
where $p$ is a holomorphic function of nonnegative real part which
does not vanish in $\Delta$.

Consider the function $$g_{1}(z):=\frac{-z}{(1-z)^{2}}\cdot
g(z)=zp(z), \quad z\in\Delta .$$ This function is the holomorphic
generator of a semigroup on $\Delta$ with its Denjoy--Wolff point
at zero.

However, the equality $$\angle\lim\limits_{z\rightarrow
1}\frac{g_{1}(z)}{z-1}=\angle\lim\limits_{z\rightarrow
1}\frac{-z}{(1-z)^{3}}\cdot g(z)=0$$ implies that $g_{1}(1)=0$ and
$g'_{1}(1)=0$. Therefore $\tau=1$, too, is the Denjoy--Wolff point
of the semigroup generated by $g_{1}$ (again by Lemma 3 in
\cite{E-S}). The contradiction we have reached proves that $g\equiv 0$ on
$\Delta$. \epr

D. M. Burns and S. G. Krantz generalized their one-dimensional
result for holomorphic self-mappings of $\Delta$ (Proposition 1)
to the open unit ball $\mathbb{B}:=\{x\in \mathbb{C}^{n}, \|x\| <1
\}$, where $\|x\|=\sqrt{|x_{1}|^{2}+|x_{2}|^{2}+...+|x_{n}|^{2}}$.

\begin{prop}[see \cite{B-K}]
Let $\mathbb{B}\subset\mathbb{C}^{n}$ be the open unit ball. Let $\Phi
:\mathbb{B}\rightarrow\mathbb{B}$ be a holomorphic mapping of the
ball to itself such that $$\Phi
(x)=\mathbf{1}+(x-\mathbf{1})+O\left(
\|x-\mathbf{1}\|^{4}\right)$$ as $x\rightarrow\mathbf{1}$. (Here
$\mathbf{1}$ denotes the distinguished boundary point
$\mathbf{1}=(1,0,...,0)$ of the ball.) Then $\Phi (x)=x$ on the
ball.
\end{prop}

At this juncture, a natural question arises: does the rigidity
result for generators (Proposition 2) admit an analogous
generalization to the open unit balls of either $\mathbb{C}^{n}$ or a
Hilbert space $H$? The following theorem gives an affirmative
answer to this question. Moreover, we show that it is sufficient
to consider the $K$-limit instead of the unrestricted one in the
assumption of the theorem.

Let $\mathbb{B}$ be the open unit ball of the Hilbert space $H$. 
For $\alpha >1$, we denote by
$$D_{\alpha}(\tau):=\left\{x\in\mathbb{B}:|1-\langle
x,\tau\rangle| <\frac{\alpha}{2}(1-\|x\|^{2})\right\}$$ the
Koranyi approach regions at $\tau\in\partial\mathbb{B}$ and say
that a mapping $f:\mathbb{B}\mapsto H$ has a $K$-limit $M$ at
$\tau$ if it tends to $M$ along every curve ending at $\tau$ and
lying in a Koranyi region $D_{\alpha}(\tau)$.

\begin{pred0}
Let $f\in\Hol (\mathbb{B},H)$ be the generator of a one-parameter
continuous semigroup on $\mathbb{B}$. If for some
$\tau\in\partial\mathbb{B}$, the $K$-limit
\begin{equation}\label{a1}
\mathrm{K}\mbox{-}\lim\limits_{x\rightarrow\tau}\frac{f(x)}{\|x-\tau\|^{3}}=0,
\end{equation}
then $f\equiv 0$ on $\mathbb{B}$.
\end{pred0}

\pr
We prove this assertion by reduction to the one dimentional case.
To this end, we fix a point $y\in\mathbb{B}$ and define the
mapping $$M_{y}(x):=\frac{y-P_{y}x-sQ_{y}x}{1-\langle x,y\rangle},
\quad x\in\mathbb{B},$$ where $P_{y}$ is the orthogonal projection
of $H$ onto the subspace generated by $y$ ($P_{0}\equiv 0$ and
$P_{y}x=\displaystyle\frac{\langle x,y\rangle}{\|y\|^{2}}y$ for
$y\neq 0$), $Q_{y}=I-P_{y}$ and $s=\sqrt{1-\|y\|^{2}}$. This
mapping is an automorphism of $\mathbb{B}$ satisfying
$M^{-1}_{y}=M_{y}$ (cf. p. 98 in \cite{G-R} and p. 34 in
\cite{R}).

Denote by $U_{y}$ a unitary operator on $\mathbb{B}$ such that
$U_{y}\tau=M_{y}\tau$. Then the mapping $m:=M_{y}\circ U_{y}$ is
an automorphism of $\mathbb{B}$ which satisfies $m(\tau)=\tau$.

Obviously, $m$ is a biholomorphism of $\mathbb{B}$ onto
$\mathbb{B}$. Therefore, by Lemma 3.7.1 on p. 30 of \cite{ERS},
the mapping
\begin{equation}\label{a2}
f_{m}(w)=[m'(w)]^{-1}f(m(w)), \quad w\in\mathbb{B},
\end{equation}
is also a holomorphic generator on $\mathbb{B}$.

Substituting
$$[m'(w)]^{-1}=[m^{-1}(x)]'_{x=m(w)}=U^{*}_{y}M'_{y}(m(w))$$ in
(\ref{a2}), we have
\begin{equation}\label{b1}
f_{m}(w)=U^{*}_{y}M'_{y}(m(w))f(m(w)), \quad w\in\mathbb{B}.
\end{equation}
Now we define a holomorphic function $g$ on the unit disk $\Delta$
of the complex plane $\mathbb{C}$ by
\begin{equation}\label{b2}
g(z):=\langle f_{m}(z\tau),\tau\rangle, \quad z\in\Delta.
\end{equation}
This function $g$ is a holomorphic generator on $\Delta$. To see
this, note that, by the Theorem in \cite{A-R-S}, the generator
$f_{m}$ satisfies the inequality $$\Re \langle
f_{m}(x)-(1-\|x\|^{2})f_{m}(0), x\rangle \geq 0 \quad \mbox{for
all} \quad x\in\mathbb{B}.$$ In particular, for $x=z\tau$, where
$z\in\Delta$, $$\Re (\langle f_{m}(z\tau),\tau\rangle
\overline{z})\geq (1-|z|^{2})\Re(\langle f_{m}(0),\tau
\rangle\overline{z}),$$  {\em i. e.}, $$\Re (g(z)\overline{z})\geq
(1-|z|^{2})\Re (g(0)\overline{z}) \quad \mbox{for all} \quad
z\in\Delta ,$$ and, consequently, by the same theorem (see
\cite{A-R-S}), $g$ is indeed a holomorphic generator on $\Delta$.
(We remark in passing that this also follows from the
characterization of generators in terms of their
$\rho$-monotonicity \cite{RS-1997}.)

We claim that under our assumptions, $g\equiv 0$ on $\Delta$.
Indeed,
\begin{equation}\label{x2}
\begin{split}g(z)&=\langle
 U^{*}_{y} M'_{y}(m(z\tau))f(m(z\tau)),\tau\rangle =\langle
 M'_{y}(m(z\tau))f(m(z\tau)),U_{y}\tau\rangle{}\\
&=\langle f(m(z\tau)),[M'_{y}(m(z\tau))]^{*}U_{y}\tau\rangle,
\quad z\in\Delta,
\end{split}
\end{equation}
and, consequently,
$$\frac{g(z)}{|z-1|^{3}}=\frac{1}{|z-1|^{3}}\langle
f(m(z\tau)),[M'_{y}(m(z\tau))]^{*}U_{y}\tau\rangle$$
$$=\frac{\|m(z\tau)-\tau\|^{3}}{|z-1|^{3}}\left\langle
\frac{f(m(z\tau))}{\|m(z\tau)-\tau\|^{3}},[M'_{y}(m(z\tau))]^{*}U_{y}\tau\right\rangle.$$

Note that each automorphism $h$ of $\mathbb{B}$ is the restriction
to $\mathbb{B}$ of a holomorphic mapping defined either on the
larger ball $B(0,R)$ centered at zero of radius $R=\displaystyle
\frac{1}{\|h^{-1}(0)\|}$ if $h(0)\neq 0$, or on all of $H$ if $h$
fixes the origin. So, $M_{y}$ and $m$ are, in fact, holomorphic
mappings defined either in the open ball $B(0,R)$ of radius
$R=\frac{1}{\|y\|}>1$ if $y\neq 0$ or on $H$ if $y=0$. Moreover,
for $z$ close enough to 1 in the nontangential approach region
$$\Gamma _{k}=\left\{z\in\Delta :\frac{|z-1|}{1-|z|}<k\right\},
\quad k>1,$$  $m(z\tau)$ belongs to the Koranyi region
$D_{\alpha}(\tau)$ whenever $\alpha >k$.

Indeed, it can be shown by direct calculations that the function
$m$ satisfies the equality $$\frac{|1-\langle
m(z\tau),\tau\rangle|^{2}}{1-\|m(z\tau)\|^{2}}=L \,
\frac{|1-z|^{2}}{1-|z|^{2}} \, , \quad z\in\Delta,$$ where
$L:=\displaystyle\frac{d}{dz}\langle m(z\tau),\tau\rangle|_{z=1}
>0$. Consequently, we have for $z\in\Gamma _{k}$, $$\frac{|1-\langle
m(z\tau),\tau\rangle|}{1-\|m(z\tau)\|^{2}}=L \,
\frac{|1-z|^{2}}{1-|z|^{2}}\cdot\frac{1}{\left|1-\langle m(z\tau
),\tau \rangle\right|}<L \, k \,
\frac{1}{1+|z|}\cdot\frac{|1-z|}{\left|1-\langle m(z\tau ),\tau
\rangle\right|} \, .$$ Since $\lim\limits_{z\rightarrow
1}\left|\frac{1-\langle m(z\tau),\tau \rangle}{1-z}\right|=L$, it
follows that if $z\in\Gamma_{k}$ is close enough to 1, then
$m(z\tau)$ is in $D_{\alpha}(\tau)$ $(\alpha >k)$. Hence, by
hypothesis (\ref{a1}) of the theorem,
$\angle\lim\limits_{z\rightarrow
1}\frac{f(m(z\tau))}{\|m(z\tau)-\tau\|^{3}}=0$.

Therefore, $\angle\lim\limits_{z\rightarrow
1}\frac{g(z)}{|z-1|^{3}}=0,$ and by Proposition 2, $g\equiv 0$ on
$\Delta$. So, by (\ref{x2}), $$\langle
f(m(z\tau)),[M'_{y}(m(z\tau))]^{*}U_{y}\tau\rangle =0 \quad
\mbox{for all} \quad z\in\Delta.$$ In particular, this equality
holds for $z=0$, {\em i. e.},
\begin{equation}\label{a5}
\langle f(y),[M'_{y}(y)]^{*}U_{y}\tau\rangle =0  \quad \mbox{for
each} \quad y\in\mathbb{B}.
\end{equation}
By direct calculations, one obtains that
\begin{equation*}
\begin{split}
&M'_{y}(x)h={}\\ &=\frac{1}{(1-\langle
x,y\rangle)^{2}}\biggl[-(1-\langle
x,y\rangle)(P_{y}+sQ_{y})h+\langle h,y\rangle
(y-P_{y}x-sQ_{y}x)\biggr].
\end{split}
\end{equation*}
Hence, $$M'_{y}(y)h=-\frac{1}{1-\|y\|^{2}}(P_{y}+sQ_{y})h,$$ and
equality (\ref{a5}) is equivalent to $$\langle
f(y),(P_{y}+sQ_{y})U_{y}\tau\rangle =0.$$ Substituting $$U_{y}\tau
=M_{y}\tau =\frac{y-P_{y}\tau -sQ_{y}\tau }{1-\langle \tau
,y\rangle}$$ in this equality, we obtain $$\langle
f(y),(P_{y}+sQ_{y})(y-P_{y}\tau-sQ_{y}\tau)\rangle =0 ,$$
 $$\langle f(y), y-P_{y}\tau -s^{2}Q_{y}\tau\rangle =0 ,$$
 $$\langle f(y), y-P_{y}\tau -(1-\|y\|^{2})(I-P_{y})\tau\rangle =0
 ,$$ and
 $$\langle f(y), y-\tau +\|y\|^{2}\tau -\langle\tau ,y\rangle y\rangle =0 \quad \mbox{for all} \quad y\in\mathbb{B}.$$
Let $y=y_{1}\tau +\tilde{y}$, where $y_{1}=\langle y,\tau\rangle$
and $\langle \tilde{y},\tau\rangle =0$.

Similarly, $f(y)=f_{1}(y)\tau +\tilde{f}(y)$ with
$f_{1}(y)=\langle f(y),\tau\rangle$ and $\langle
\tilde{f}(y),\tau\rangle =0$ for all $y\in\mathbb{B}$.

Using this notation, we have  $$\langle f_{1}(y)\tau , y_{1}\tau
-\tau +\|y\|^{2}\tau -|y_{1}|^{2}\tau \rangle =-\langle
\tilde{f}(y),\tilde{y}-\overline{y}_{1}\tilde{y}\rangle$$ and
$$(1-\overline{y}_{1}-\|\tilde{y}\|^{2})f_{1}(y)=(1-y_{1})\langle
\tilde{f}(y),\tilde{y} \rangle.$$ Differentiating this equality
with respect to $\overline{y}_{1}$, we conclude that it can hold
only if $f_{1}(y)=0$ and
\begin{equation}\label{a6}
\langle \tilde{f}(y), \tilde{y}\rangle =0 \quad \mbox{for
all} \quad y\in\mathbb{B}.
\end{equation}

Now let $\sigma$ be an arbitrary unit vector orthogonal to $\tau$,
{\em i. e.}, $\langle \sigma ,\tau\rangle =0$. Suppose that
$\tilde{y}=y_{2}\sigma +u$, where $y_{2}=\langle
\tilde{y},\sigma\rangle$ and $\langle u,\sigma\rangle =0$.

Similarly, $\tilde{f}(y)=f_{2}(y)\sigma +v(y)$ with
$f_{2}(y)=\langle \tilde{f}(y),\sigma \rangle$ and $\langle
v(y),\sigma\rangle =0$ for all $y\in\mathbb{B}$. Then by
(\ref{a6}), $$f_{2}(y)\overline{y}_{2}=-\langle v(y),u\rangle.$$
Differentiating this equality with respect to $\overline{y}_{2}$,
we obtain $f_{2}(y)=0$. Hence, $f\equiv 0$ on $\mathbb{B}$. \epr

Following L. A. Harris \cite{Har}, we define the numerical range
of each $h\in\Hol (\mathbb{B} , H)$ which has a norm continuous
extension to $\overline{\mathbb{B}}$ by $$V(h):=\{\langle
h(x),x\rangle :\| x \| =1\}.$$

For an arbitrary holomorphic mapping $h\in \Hol (\mathbb{B},H)$
and for each $s~\in ~(0,1)$, we define the mapping
$h_{s}:\frac{1}{s}\mathbb{B}\mapsto H$ by $$h_{s}:=h(sx), \quad
\|x\|<\frac{1}{s},$$ and set $$L(h):=\lim\limits_{s\rightarrow
1^{-}}\sup\Re(V(h_{s})).$$

It is known (Theorem 1 in \cite{H-R-S}) that the mapping $I-h$ is
a generator if and only if $L(h)\leq 1$. So the following
corollary is an immediate consequence of our theorem.

\begin{corol}
Let $h\in\Hol (\mathbb{B},H)$ with $L(h)\leq 1$. If for some
$\tau\in\partial\mathbb{B}$, the $K$-limit
\begin{equation}\label{a10}
\mathrm{K}\mbox{-}\lim\limits_{x\rightarrow\tau}\frac{h(x)-x}{\|x-\tau\|^{3}}=0,
\end{equation}
then $h\equiv I$ on $\mathbb{B}$.
\end{corol}

Since obviously $L(h)\leq 1$ for all self-mappings of $\mathbb{B}$,
this corollary properly contains Proposition 3.

{\bf Acknowledgment.} The third author was partially supported by
the Fund for the Promotion of Research at the Technion and by the
Technion President's Research Fund.

\end{document}